\pdfoutput=1
\documentclass{amsart}

\usepackage[pdftex]{graphicx}
\usepackage[english]{babel}
\usepackage{amssymb}
\usepackage{amsfonts}
\usepackage{amsmath}
\usepackage{float}
\usepackage{url}
\usepackage{mathrsfs}
\usepackage{parskip}

\graphicspath{ {./Figures/} }

\theoremstyle{definition}

\theoremstyle{remark}

\numberwithin{equation}{section}

\begin{document}

\title{Fatigue Estimation Methods Comparison \\ for Wind Turbine Control}

\author{J.J.~Barradas~Berglind}
%\address{Automation \& Control, Department of Electronic Systems, Aalborg
%University, Aalborg East 9220, DK.}
%\curraddr{Department of Mathematics and Statistics,
%Case Western Reserve University, Cleveland, Ohio 43403}
%\email{jjb@es.aau.dk}
\thanks{}

\author{Rafael~Wisniewski}
%\address{Automation \& Control, Department of Electronic Systems, Aalborg
%University, Aalborg East 9220, DK.}
%\email{raf@es.aau.dk}
\thanks{\emph{Pre-print submitted to Wind Energy.}}

%    General info
%\subjclass[2000]{Primary 54C40, 14E20; Secondary 46E25, 20C20}

%\date{January 1, 2001 and, in revised form, June 22, 2001.}

\dedicatory{Automation \& Control, Department of Electronic Systems \\ Aalborg
University, Aalborg East 9220, DK. e-mail: \{jjb,raf\}@es.aau.dk}

\keywords{Fatigue Estimation, Load Estimation, Fatigue for Wind Turbine Control, Rainflow Counting, Spectral Methods, Hysteresis Operator}

\begin{abstract}
Fatigue is a critical factor in structures as wind turbines exposed to harsh operating conditions, both in the design stage and control during their operation. In the present paper the most recognized approaches to estimate the damage caused by fatigue are discussed and compared, with special focus on their applicability for wind turbine control. The aim of this paper is to serve as a guide among the vast literature on fatigue and shed some light on the underlying relationships between these methods.
\end{abstract}

\maketitle

%\section*{This is an unnumbered first-level section head}
%This is an example of an unnumbered first-level heading.

%% The correct journal style for \specialsection is all uppercase; a known bug
%% in amsart.cls prevents this, so input must be uppercase until it is fixed.
%\specialsection*{This is a Special Section Head}
%\specialsection*{THIS IS A SPECIAL SECTION HEAD}
%This is an example of a special section head%
%%%%%%%%%%%%%%%%%%%%%%%%%%%%%%%%%%%%%%%%%%%%%%%%%%%%%%%%%%%%%%%%%%%%%%%%
%\footnote{Here is an example of a footnote. Notice that this footnote
%text is running on so that it can stand as an example of how a footnote
%with separate paragraphs should be written.
%\par
%And here is the beginning of the second paragraph.}%
%%%%%%%%%%%%%%%%%%%%%%%%%%%%%%%%%%%%%%%%%%%%%%%%%%%%%%%%%%%%%%%%%%%%%%%%

\section{Introduction and Motivation} 

Fatigue has been widely and exhaustively studied from different perspectives, and the literature is vast and approached from different perspectives; thus, incorporating fatigue or wear in components of a wind turbine in a control problem may seem as a daunting task. Fatigue is regarded as a critical factor in structures such as wind turbines, where it is necessary to ensure a certain life span under normal operating conditions in a turbulent environment. These environmental conditions lead to irregular loadings, which is also the case for waves and uneven roads. The main focus of the present is on fatigue estimation methods for wind turbine control, and as such the most widely used methods are described, with special emphasis in the applicability of these techniques for control.

In general, fatigue can be understood as the weakening or breakdown of a material subject to stress, especially a repeated series of stresses. From a materials perspective, it can be also thought of as elastoplastic deformations causing damage on a certain material or structure, compromising its integrity. 

Fatigue is a phenomenon that occurs in a microscopic scale, manifesting itself as deterioration or damage. Consequently, it has been of interest in different fields and has been studied extensively with different perspectives; a very detailed history of fatigue can be found in \cite{schutz96}. It could be argued that two major turning points in the history of fatigue came firstly with the contributions of W{\"o}hler, who suggested design for finite fatigue life in the 1860's \cite{Wohler_1860} and the so-called W{\"o}hler curve (or S-N curve stress versus number of cycles to failure) which still sets the basis for theoretical damage estimation; and secondly with the linear damage accumulation rule by Palmgren \cite{Palmgren_1924} and Miner \cite{Miner_1945}, still under use nowadays. 

\section{Fatigue Estimation for Wind Turbine Control}

Perhaps the most recognized and used measure for fatigue damage estimation is the so-called rainflow counting (RFC) method, which is used in combination with the Palmgren-Miner rule. In the wind turbine context, the impact on fatigue from a load can be described by an equivalent damage load (EDL); basically, the EDL is calculated using the Palmgren-Miner rule to determine a single, constant-rate fatigue load that will produce equivalent damage \cite{Sutherland99}.

Load or fatigue reduction techniques for wind turbines can be roughly divided in active and passive. The former makes use of the controller, e.g., by changing the pitching angle or the generator torque, while the latter entails the design of the structure. In \cite{bottasso13}, both strategies are combined to reduce loads in the blades. In the wind turbine control context, the control algorithm may have substantial effects on the wind turbine components; for example, controlling the pitching angle may lead to thrust load changes, which consequently affects the loads on the tower and blades \cite{bossanyi03}. In \cite{bossanyi03}, \cite{bossanyi05}, \cite{larsen05}  reductions in loading are achieved by controlling the pitch of each blade independently; the damage of different control strategies is assessed by EDL, using S-N curves. In \cite{lescher06}, \cite{nourdine10} a load reduction control strategies are proposed, where the damage is evaluated using the RFC algorithm. Model predictive control (MPC) strategies using wind preview have been proposed in \cite{Sol_11}, \cite{madsen12} to reduce loads, evaluated via EDL. In \cite{Hamm_07}, control strategies were designed, by approximating fatigue load by an analytical function based on spectral moments. The Aeolus project \cite{aeolus} has a simulation platform, which considers the fatigue load of wind farm for optimization as a post-processing method. 

A large amount of the current control methods rely on the calculation of the damage either by EDL or RFC, which can be only used as post-processing tools; other methods are based on minimization of some norms of the stress on different components of the wind turbine, which are hoped to reduce fatigue, but they are not a reliable characterization of the damage \cite{Sol_11}, \cite{Mirzaei_13}. Thus, in this paper we will introduce and compare the most recognized fatigue estimation methods, and explore different alternatives with a focus on whether they can be incorporated in control loops and thus be used in the controller synthesis directly.

\section{Fatigue Estimation Methods}

Some of the most recognized approaches to estimate the damage caused by fatigue will be discussed and compared in the sequel. From a materials perspective, an extensive survey for homogeneous materials was done in \cite{fatemi_98}. In the wind turbine context, \cite{Sutherland99} goes through the counting and spectral techniques used for wind turbine design. The perspective taken here is from a control point of view and as such we categorize the fatigue estimation methods as follows:
\begin{enumerate}
	\item{Counting methods} 
	\item{Frequency domain or spectral methods} 
	\item{Stochastic methods} 
	\item{Hysteresis operator} 
\end{enumerate}

In all cases, we assume that the input signal is obtained from time history of the loading parameter of interest, such as force, torque, stress, strain, acceleration, or deflection \cite{fatiguelee_2005}. 

\subsection{Counting Methods}

Cycle counting methods are algorithms that identify fatigue cycles by combining and extrapolating information from extrema (maxima and minima) in a time series. These algorithms are used together with damage accumulation rules, which calculate the total damage as a summation of increments. The most popular method among the counting methods is the so-called rainflow counting (RFC) method, jointly with the Palmgren-Miner rule of linear damage accumulation to calculate the expected damage. The Palmgren-Miner rule is the most popular due to its simplicity; however, by applying it one assumes a fixed-load, neglecting interaction and sequence effects that might have a significant contribution to the damage, e.g., \cite{Agerskov} for tests with random loading.

Other cycle counting methods include: peak-valley counting (PVC), level-crossing counting (LCC), range counting (RC), and range-pairs counting (RPC); for more details see \cite{ASTMe1049} and \cite{Benasciutti_Thesis}. Here, we will focus on the RFC method, which is the most widely used and the  most accurate in identifying the damaging effects caused by complex loadings, \cite{Dowling71}. The rainflow counting method, first introduced by Endo \cite{Endo_1967}, has a complex sequential and nonlinear structure in order to decompose arbitrary sequences of loads into cycles, and its name comes from an analogy with roofs collecting rainwater to explain the algorithm, sometimes also referred to as pagoda roof. A figure depicting the described procedure is shown below in Figure \ref{fig:Int_RFC}.

\begin{figure}[htb]
\centering
\includegraphics[scale=1.1]{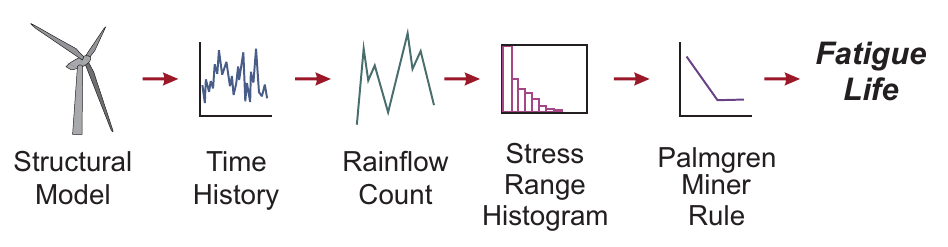}
\caption{Rainflow counting damage estimation procedure.}
\label{fig:Int_RFC}
\end{figure}

For many materials there is an explicit relation between number of cycles to failure and cycle amplitude, which is known as S-N or W{\"o}hler curves, given as a line in a log-log scale as

\begin{align}
	s^{k}N = K,
\label{SNcurve}
\end{align}

where $k$ and $K$ are material specific parameters and $N$ is the number of cycles to failure at a given stress amplitude $s$. Then, for a time history, the total damage under the linear accumulation damage (Palmgren-Miner) rule is given as

\begin{align}
	D(T) = \sum\limits_{i=1}^{N(T)}\Delta D_{i} = \sum\limits_{i=1}^{N(T)}\frac{1}{N_{i}},
\label{DamagePM1}
\end{align}

for damage increments $\Delta D_{i}$ associated to each counted cycle, $N_{i}$ the number of cycles to failure associated to stress amplitude $s_{i}$, and the number of all counted cycles $N(T)$. Taking the S-N curve relationship in \eqref{SNcurve}, we can rewrite \eqref{DamagePM1} as

\begin{align}
	D(T) = \sum\limits_{i=1}^{N(T)}\frac{s_{i}^{k}}{K}.
\end{align}

\begin{figure}[htb]
\centering
\includegraphics[width=0.95\textwidth]{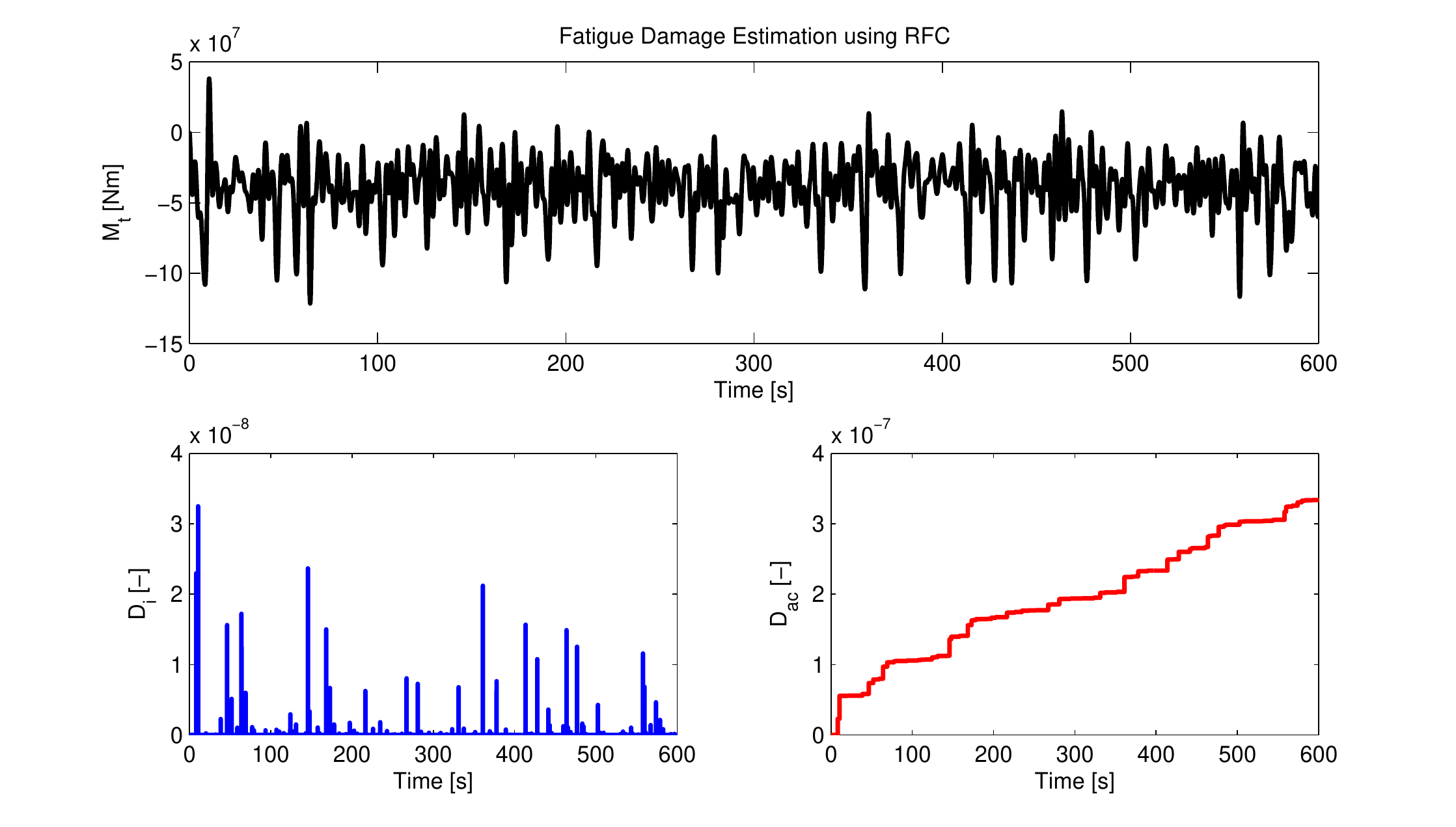}
\caption{Rainflow counting algorithm example, using the toolbox from \cite{Nieslony_09}.}
\label{fig:Ex1_RFC}
\end{figure}

Different RFC algorithms have been proposed such as \cite{downing1982} and \cite{rychlik_3}, with different rules but providing the same results. A way to implement the RFC algorithm is using the Rainflow toolbox introduced in \cite{Nieslony_09}. An example is presented below, using the wind turbine model from the standard NREL 5MW wind turbine \cite{Jonk_NREL5MW}, running is closed-loop with standard pitch and torque controllers. The input used for the comparison is a time series of the tower bending moment extracted after the simulation of $600$ seconds. The results are presented on Figure \ref{fig:Ex1_RFC}. On the top the input stress is shown, and in the bottom part the instantaneous damage and the accumulated damage are shown. For our example, we will let $k=4$ and $K=6.25\times 10^{37}$ as in \cite{Hamm06_Thesis}, where the value of $k$ is adequate for steel structures. For this example, the instantaneous damage was extrapolated to its causing time, such that it can be plotted in the right time scale instead of the  reduced turning-point scale.

\begin{figure}[htb]
\centering
\includegraphics[scale=0.7]{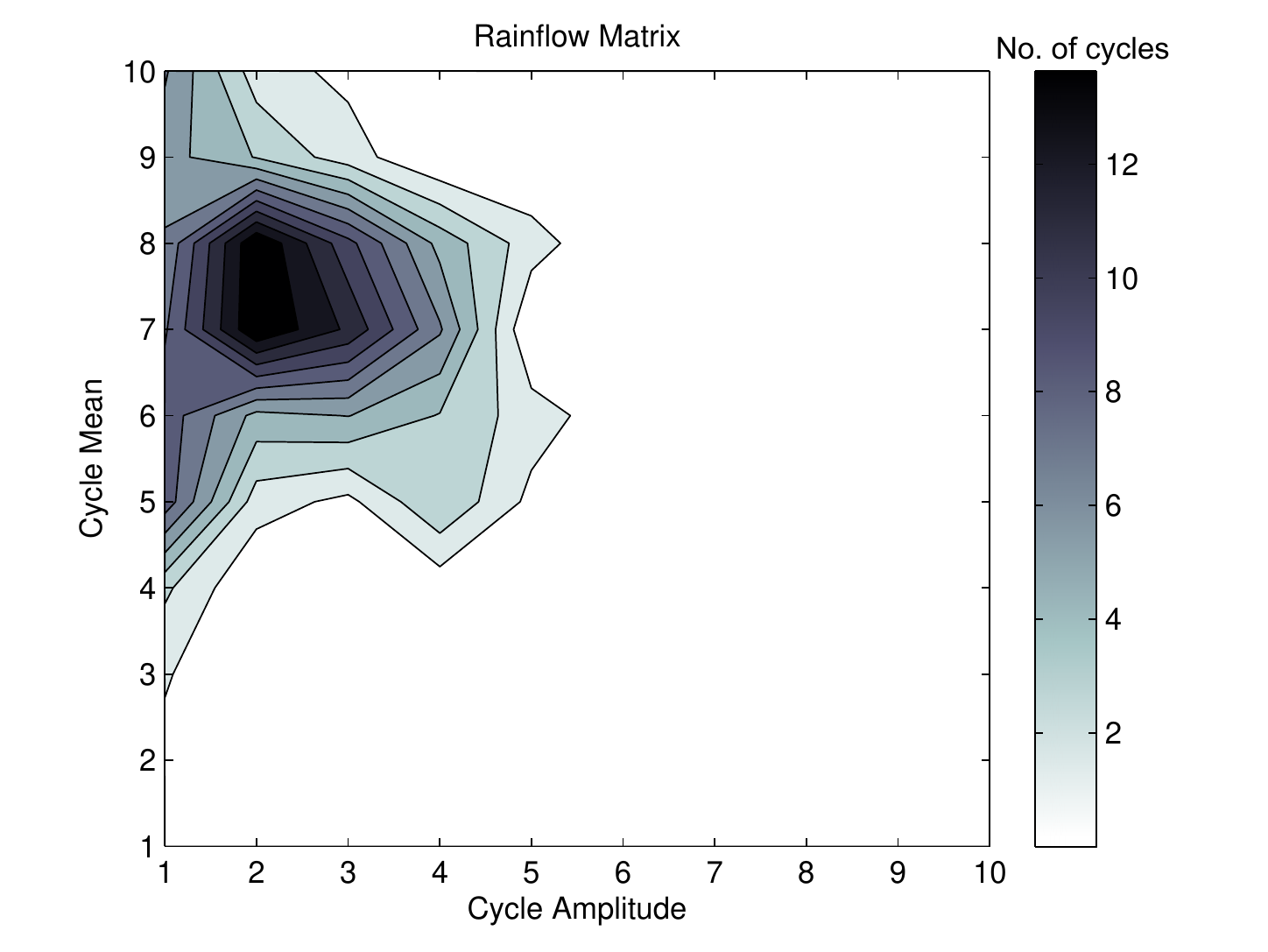}
\caption{Rainflow Matrix, using the toolbox from \cite{Nieslony_09}.}
\label{fig:Ex1_RFC2}
\end{figure}

Other outputs provided by the toolbox in \cite{Nieslony_09} are amplitude and cycle mean histograms, as well as the so-called rainflow matrix (RFM), from which the number of counted cycles with a given amplitude and mean value are obtained from the given stress history. Since the RFM will play a role further on this paper, we will elaborate on its construction. Load signals can be discretized to a certain number of levels, allowing an efficient storage of the cycles in a so-called rainflow matrix, which is an upper triangular matrix by definition. Consequently, cycle amplitudes and mean values can be grouped in bins, such that the cycle count can be summarized as a matrix (for details see \cite{rychlik_2}, and Chapter 2 in \cite{Benasciutti_Thesis}); sometimes this matrix is shown transposed. The rainflow matrix for the aforementioned example is depicted on Figure \ref{fig:Ex1_RFC2} for 10 bins, where cycle mean is on the $y-$axis, cycle amplitude in the $x-$axis and number of cycles in the $z-$axis.

Lastly, NREL has a an estimator of fatigue-life called \texttt{MLife} (currently in alpha version, an improvement on \texttt{MCrunch} \cite{MCrunch}), which runs the RFC algorithm of \cite{Nieslony_09}. \texttt{MLife} calculates fatigue life for one or several time series, incorporating the Goodman correction to the damage calculation (to account and correct for the fixed-load assumption). These calculations include short-term damage equivalent loads and damage rates, lifetime results based on time series, accumulated lifetime damage, and time until failure \cite{MLife}.

\subsection{Spectral Methods}

\begin{figure}[htb]
\centering
\includegraphics[scale=1.1]{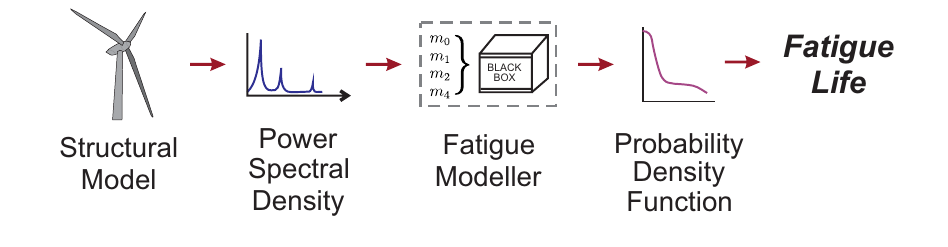}
\caption{Spectral methods damage estimation procedure.}
\label{fig:Int_Spectral}
\end{figure}

An alternative to counting methods are the so-called spectral or frequency domain methods \cite{Bishop_99}, which assume narrow band processes and calculate the lifetime estimate by using an empirical formula that uses the spectral moments of the input signal; the aim of these methods is to approximate the rainflow density of the RFC algorithm. This procedure is depicted on Figure \ref{fig:Int_Spectral}. It is worth mentioning that some of these methods are based on empiric formulas, being essentially black-box and may be restricted to Gaussian histories. A comparison of different spectral methods was carried out in \cite{Benasciutti_06}.

Spectral methods are based on statistical information of the signal of interest, i.e., its spectral moments. Following from \cite{Bishop_99} and \cite{Hamm_07}, the \emph{m}$^{th}$ spectral moment of the process $x(t)$ is defined as

\begin{align}
	\lambda^{x}_{m} = \frac{1}{\pi}\int\limits_{0}^{\infty}f^{m}\cdot S_{x}(f) df,
\end{align}

where $S_{x}(f)$ is the power density (PSD) of the process, with the following properties

\begin{align}
	\lambda^{x}_{0} = \sigma^{2}_{x}, \;\; \lambda^{x}_{2} = \sigma^{2}_{\dot{x}} \;\; \text{and} \;\; \lambda^{x}_{4} = \sigma^{2}_{\ddot{x}}.
	\label{moments4}
\end{align}

In other words, the variance of the process is given by $\lambda^{x}_{0}$, the variance of the process' first derivative is then given by the second moment, and lastly the variance of the process' second derivative is given by the fourth moment. Consequently, following the results in \cite{Benasciutti05} and \cite{Ryc93} the damage rate for narrow-banded Gaussian stress histories is given by

\begin{align}
d_{\curlywedge} = \frac{1}{2\pi}\sqrt{\frac{\lambda_{4}}{\lambda_{2}}}\frac{1}{K}\left(2\sqrt{2\lambda_{0}}\right)^{k}\Gamma\left(1+\frac{k}{2}\right),
\end{align}

where $\Gamma(\cdot)$ corresponds to the gamma distribution, and $k$, $K$ are the S-N parameters used in the RFC case. In \cite{Benasciutti05}, the authors proposed an estimate of the expected fatigue damage rate given as the narrow-band approximation augmented with a correction factor to account for the process not necessarily being narrow-band

\begin{align}
E\left[d\right] \approx  d_{\curlywedge} \cdot \left(b+(1-b)\alpha_{2}^{k+1}\right) 
\end{align}

with 

\begin{align}
	b = \frac{\left(\alpha_{1}-\alpha_{2}\right)\left[1.112\left(1+\alpha_{1}\alpha_{2}-(\alpha_{1}+\alpha_{2})\right)e^{2.11\alpha_{2}}+\left(\alpha_{1}-\alpha_{2}\right)\right]}{\left(\alpha_{2}-1\right)^{2}}
\end{align}

and 

\begin{align}
	\alpha_{1} = \frac{\lambda_{1}}{\sqrt{\lambda_{0}\lambda_{2}}}, \hspace{12pt} \alpha_{2}= \frac{\lambda_{2}}{\sqrt{\lambda_{0}\lambda_{4}}}.
\end{align}

In \cite{Hamm_07} and \cite{Hamm06_Thesis}, the numerical integration of the spectral density as in \eqref{moments4} is avoided, since the spectral moments are computed by means of polynomial evaluation and differentiation, involving a logarithm and an inverse tangent function. This allowed the method to be incorporated in the control loop.

In order to compare the spectral method with the example presented in the previous section, the spectral moments $\lambda=(\lambda_{0},\lambda_{1},\lambda_{2},\lambda_{4})$ of the time series were calculated using the \texttt{WAFO toolbox} \cite{WAFO} (through integration)

\begin{align}
	\lambda=\{4.4071E^{14},-3.949E^{07},2.2904E^{11},2.1263E^{11}\},
\end{align}

and then the damage was computed using the Benasciutti approximation, using the \texttt{Matlab} script in Appendix B.3. of \cite{Hamm06_Thesis}, such that

\begin{align}
	d_{B} = 4.0024E^{-12},
\end{align}

which is a little off compared to the RFC case; this can be explained by the fact that we need to scale the damage rate according to the geometry of the system, which is generally unknown. However, the obtained damage rate can be normalized to be used for control purposes, for details see \cite{Hamm_07}. In \cite{Ragan_07} the RFC method is compared with the spectral method using Dirlik's formula (which approximates the rainflow density, see \cite{Dirlik_85}) for fatigue analysis of several components of wind turbines, where it is concluded that spectral methods work very well in some cases, but rather poorly in others due to the narrow band assumption. However, spectral methods do have the advantage of conveniently relying on spectral information that is easier to estimate from limited data. 

\subsection{Stochastic Methods}

In \cite{Sobczyk_87}, a thorough survey of stochastic methods for fatigue estimation in materials is presented, including reliability-inspired approaches, evolutionary probabilistic approaches and models for random fatigue crack growth. Modeling fatigue as a stochastic process makes sense due to the random nature of fatigue, which becomes more obvious under time-varying random loading. 

Due to the broadness of this class of methods, we will focus on one example of the evolutionary approach. Following \cite{Sobczyk_87}, by introducing the hypothesis that the process is Markovian, such that future outcomes only depend on present information, disregarding the past. This way, we will have a random process with only forward transitions,

\begin{align}
E_{0} \rightarrow E_{1} \rightarrow \cdots \rightarrow E_{k} \rightarrow E_{k+1} \cdots \rightarrow E_{n}=E^{*},
\label{Eq:fatigueMC}
\end{align}

where $E_{0}$ denotes a damage-free state and $E^{*}$ characterizes the ultimate damage or destruction. Letting $P_{k}(t)$ be the probability that the specimen at time $t$ is on state $E_{k}
$ (notice that the state transitions are discrete, while the time evolution is continuous), then we obtain the following system of differential equations

\begin{align}
\frac{dP_{0}(t)}{dt} &= q_{0}P_{0}(t)  \nonumber\\
\frac{dP_{k}(t)}{dt} &= q_{k}P_{k}(t) + q_{k-1}P_{k-1}(t), \hspace{12pt} k \geq 1,
\end{align}

or in shorter notation

\begin{align}
\frac{dP_{k}(t)}{dt} &= Q P_{k}(t), \hspace{12pt} k \geq 0,
\label{Mchain}
\end{align}

which corresponds to a Markov chain (MC) with intensity or transition matrix $Q$. Markov chains are well studied and have been successfully used in control settings; however, a shortcoming of this approach is that it is assumed that the intensity matrix $Q$ is not generally known. It could be assumed that the intensities are obtained from physical experiments, but this would correspond to a certain load; so, if the load changes, the parameters will change as well. However, the elements of $Q$ could be identified, using for instance recursive maximum likelihood identification methods, in order to capture the shifts in the load introduced by the controller.

\begin{figure}[htb]
\centering
\includegraphics[scale=0.6]{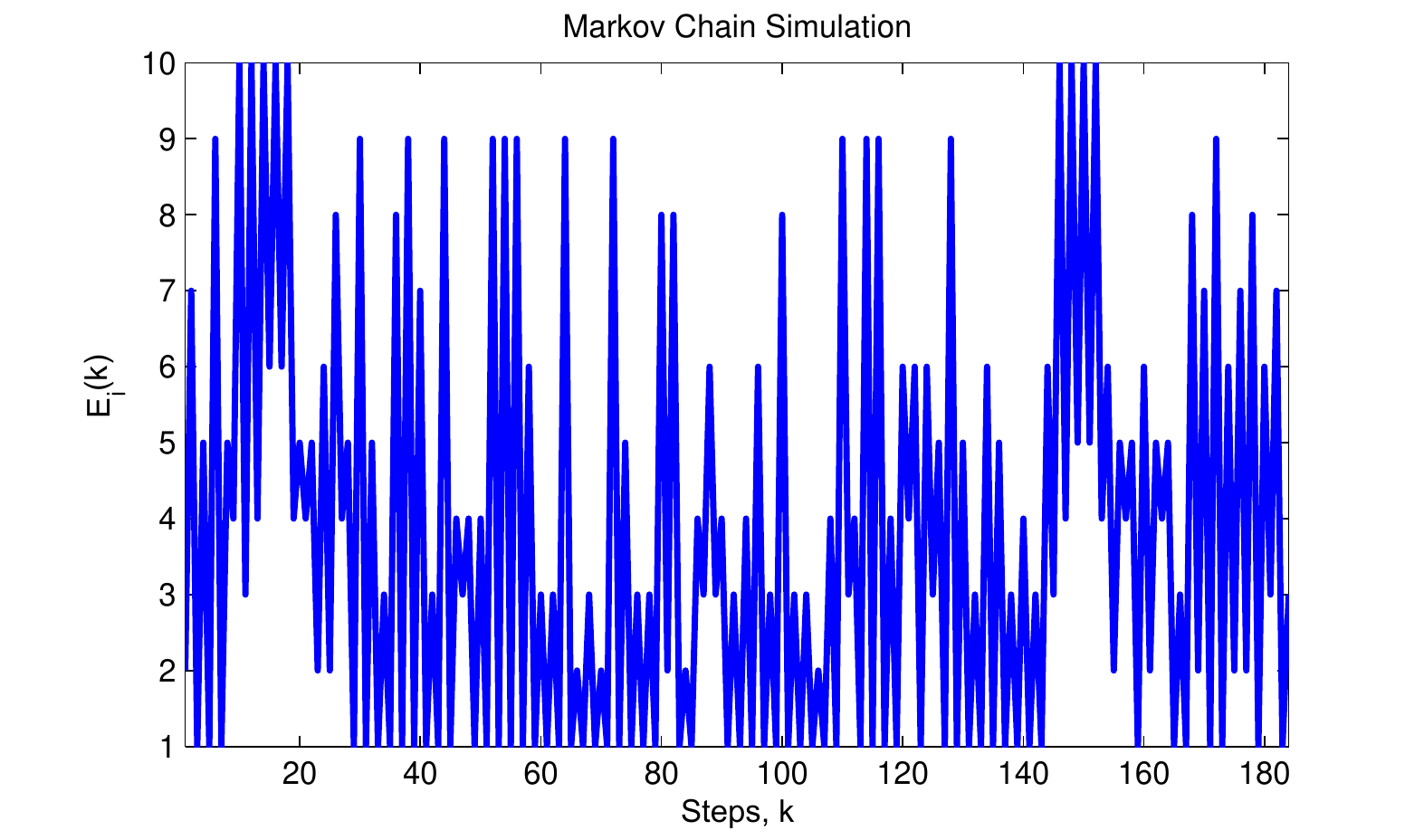}
\caption{Markov Chain simulation, using the WAFO toolbox \cite{WAFO}.}
\label{fig:MCsimulation}
\end{figure}

In the present, for the sake of comparison, we will make use of the equivalence in \cite{rychlik_2}, where a method to convert between rainflow matrix to a Markov matrix is presented. As an example, we take the rainflow matrix depicted in Figure \ref{fig:Ex1_RFC2} and use the \texttt{WAFO} toolbox to convert it into a Markov matrix, and obtain its corresponding intensity matrix $Q$. Additionally, the MC is simulated for as many steps as the length of turning points of the RFC algorithm, such that the instantaneous damage can be reconstructed in the appropriate time instances. The simulation of the MC is presented on Figure \ref{fig:MCsimulation}, where the size of the MC corresponds to the number of bins of the RFM.

\begin{figure}[htb]
\centering
\includegraphics[scale=0.8]{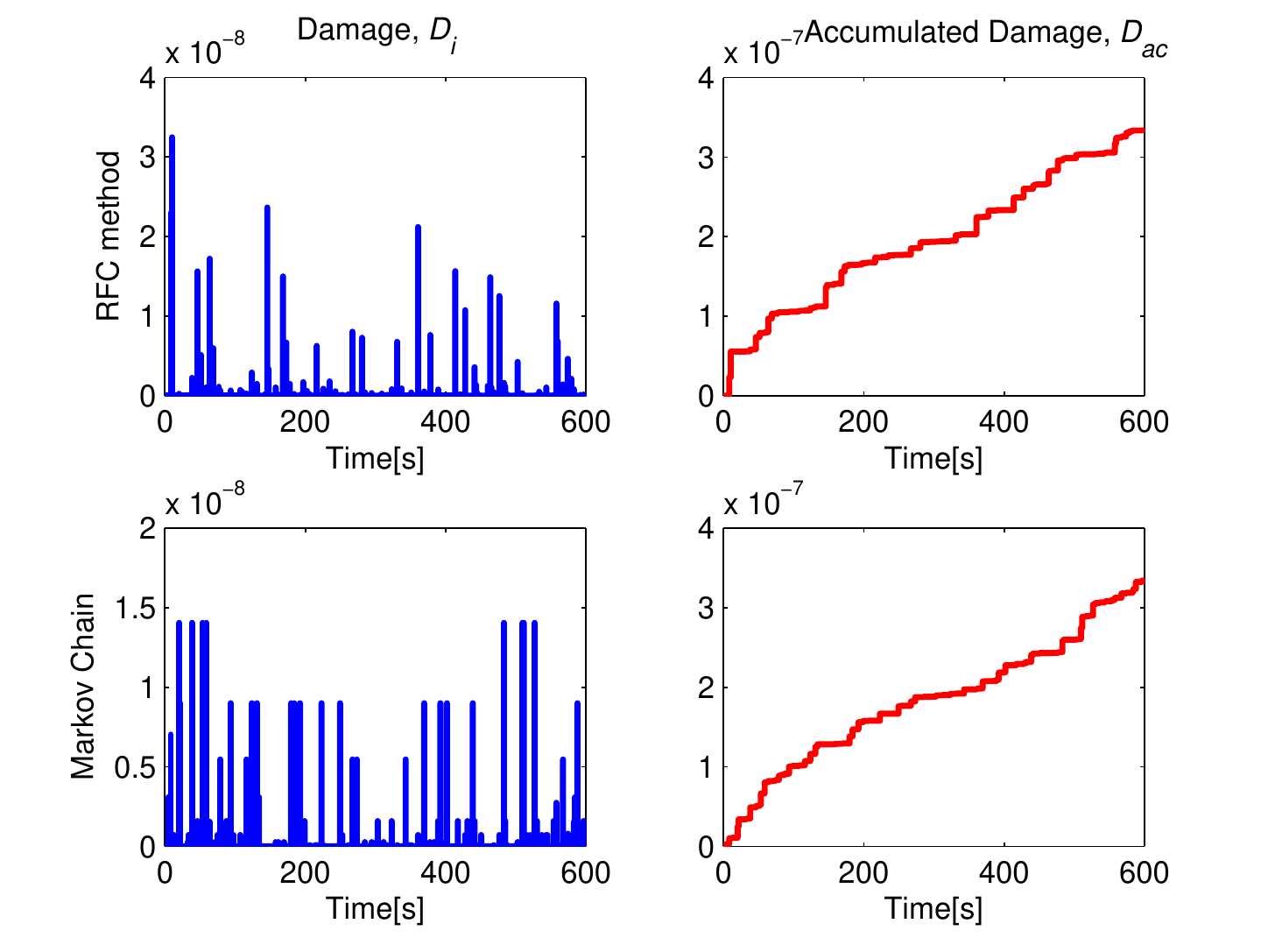}
\caption{RFC versus Markov chain method damage comparison.}
\label{fig:RFCvsMC}
\end{figure}

Then, the damage evolution is scaled according to the RFM amplitudes, and afterwards the Palmgren-Miner rule is used. One of the possible realizations is compared against the RFC method on Figure \ref{fig:RFCvsMC}. Note that many realizations for the damage evolution are possible, since the MC in \eqref{Mchain} is governed by probabilities.

\subsection{Hysteresis Operator}

As mentioned in \cite{downing1982} and \cite{fatiguelee_2005}, the purpose of the RFC method is to identify the closed hysteresis loops in the stress and strain signals. In \cite{tchankov1998}, an incremental method for the calculation of dissipated energy under random loading is presented, where the dissipated hysteresis energy to failure is used as the fatigue life parameter; the physical interpretation is that as some of the energy is dissipated, certain damage is introduced to a material or structure. 

In \cite{BroSpre_96} an equivalence between symmetric RFC and a Preisach hysteresis operator is provided. This is a very useful result, since it gives the opportunity to incorporate the fatigue estimation online in the control loop. Additionally, this method is strongly related to the physical behavior of the damaging process as explained in \cite{BroDreKre_96}. If one associates values to individual cycles or hysteresis loops, it is being assumed that the underlying process is rate independent, thus meaning that only the loops themselves are important, but not the speed with which they are traversed; in other words,  what causes the damage is the cycle amplitude and not how fast it occurs. Rate independent processes are mathematically formalized as hysteresis operators, see \cite{KrasPok_89}, \cite{Mayergoyz_91} \cite{BroSpre_96}.

The aforementioned equivalence in \cite{BroSpre_96} between symmetric rainflow counting (RFC) and a type of Preisach operator, is given as

\begin{align}
		D_{ac}(s)=\sum_{\mu<\tau}\frac{c(s)(\mu,\tau)}{N(\mu,\tau)}=\text{Var}(\mathcal{W}(s)).
\label{Cor2_13}
\end{align}

where the left-hand side corresponds to the damage given by the RFC with $c(s)(\mu,\tau)$ being the rainflow count associated with a fixed string $s=(v_{0},\cdots,v_{N})$, counting between the values of $\mu$ and $\tau$, and $N(\mu,\tau)$ denotes the number of times a repetition of the input cycle $(\mu,\tau)$ leads to failure. 

The right-hand side of \eqref{Cor2_13} is the variation of a special hysteresis operator, namely the Preisach operator defined as, 

\begin{align}
		\mathcal{W}(s) = \int_{\mu<\tau}\rho(\mu,\tau) \mathcal{R}_{\mu,\tau}(s)d\mu d\tau.
		\label{preisach_op}
\end{align}

with density function $\rho(\mu,\tau)$, interpreted as a gain that changes with the different values of $\mu$ and $\tau$, being a function of $N(\mu,\tau)$. To interpret the right-hand side of \eqref{Cor2_13} we will need to introduce the relay operator $\mathcal{R}_{\mu,\tau}(s) = \mathcal{R}_{\mu,\tau}(v_{0},\cdots,v_{N})=(w_{0},\cdots,w_{N})$, 
where its output is given by

\begin{align}
		w_{i} = \left\{
  \begin{array}{l l}
    1, & \quad v_{i}\geq \tau,\\ 
		0, & \quad v_{i}\leq \mu,\\ 
    w_{i-1}, & \quad \mu<v_{i}<\tau.
  \end{array} \right.
\end{align}

with $\mu<\tau$ and $w_{-1} \in \{0,1\}$ given. The relevant threshold values for the relays $\mathcal{R}_{\mu,\tau}$ in the Preisach operator $\mathcal{W}(s)$ then lie within the triangle

\begin{align}
		P=\left\{(\mu,\tau)\in\mathbb{R}^{2},-M\leq \mu\leq \tau\leq M\right\}.
\label{prei_plane}
\end{align}

known as the Preisach plane. The variation operator $\text{Var}(\cdot)$ is a counting element defined as

\begin{align}
		\text{Var}(s)=\sum^{N-1}_{i=0}\left|v_{i+1}-v_{i}\right|
\end{align}

for an arbitrary input sequence $s=(v_{0},\cdots,v_{N})$; so essentially, $\text{Var}(\mathcal{W}(s))$ represents the counting between the thresholds $\mu$ and $\tau$, weighted by certain gain $\rho$. Notice as well, that the limit under the integral defining the Preisach operator in \eqref{preisach_op} is congruent with the RFM being upper triangular. 

In order to apply this fatigue estimation method to the previous example, the Preisach operator $\mathcal{W}(s)$ was approximated as a parallel connection of three relay operators

\begin{align}
\mathcal{H}(s)=\sum_{i}\nu(\mu_{i},\tau_{i})\mathcal{R}_{\mu_{i},\tau_{i}}(s),
\end{align}

for $i=\{1,2,3\}$. The thresholds were set to $(\mu_{1},\tau_{1})=(-0.66M,0.66M)$, $(\mu_{2},\tau_{2})=(0.66M,0.66M)$ and $(\mu_{3},\tau_{3})=(-0.66M,-0.66M)$ corresponding to uniform discretization, where $M$ is the bound for the Preisach plane in \eqref{prei_plane} calculated as $M=\max\left\{\min\left\{s\right\},\max\left\{s\right\}\right\}$. The initial conditions of the relays were given according to the following condition:

\begin{align}
		w_{-1}(\mu_{i},\tau_{i}) = \left\{
  \begin{array}{l l}
    1, & \quad \mu_{i}+\tau_{i}<0, \\ 
		0, & \quad \mu_{i}+\tau_{i}\geq 0.
  \end{array} \right.
\end{align}

Lastly, since the Preisach density function $\rho(\mu,\tau)$, captured by the weightings on each relay $\nu(\mu_{i},\tau_{i})$ is unknown, the individual weightings of each relay were normalized such that $\nu_{1}=\alpha$, $\nu_{2}=\alpha^2$, $\nu_{3}=\alpha^3$ for $\nu_{1}+\nu_{2}+\nu_{3}=1$. Thus the accumulated damage can be written in closed form as

\begin{align}
	D_{ac}(s) = \text{Var}\left(\mathcal{H}(s)\right),
\label{eq:EXdamage}
\end{align}

where we let the input signal $s$ be the tower bending moment from the previous examples. %, such that $s\equiv M_{t}$.

\begin{figure}[htb]
\centering
\includegraphics[scale=0.65]{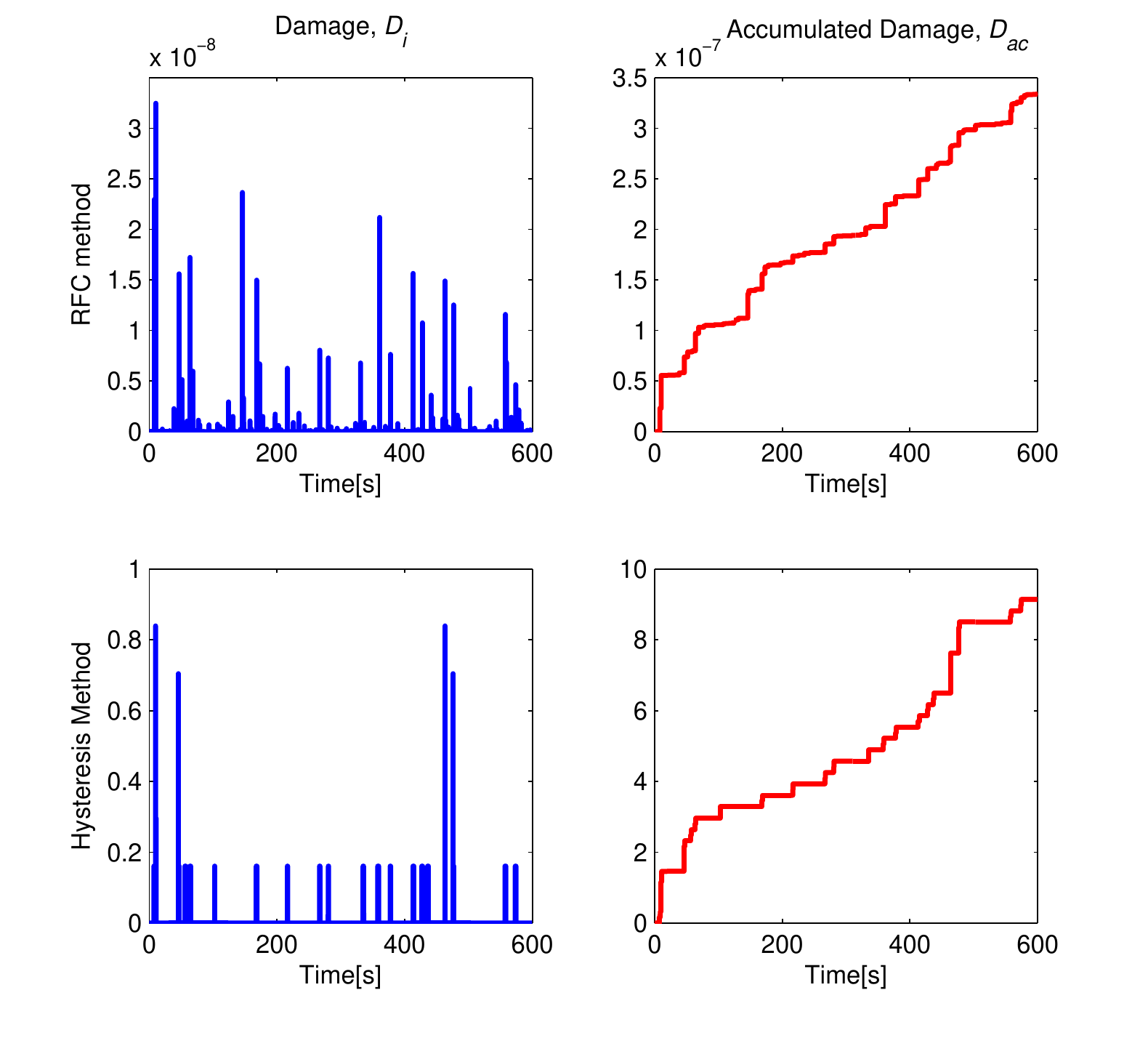}
\caption{RFC versus Hysteresis method damage comparison.}
\label{fig:RFCvsHyst}
\end{figure}

A comparison between the RFC, using the procedure described before, and the hysteresis method obtained by \eqref{eq:EXdamage} is shown in Figure \ref{fig:RFCvsHyst}. Even though the magnitude in the damage given by the hysteresis method is off scale, this could be resolved by identifying the Preisach density, see \cite{Kris_01} for an identification procedure and a summary of other identification methods.

It is worth mentioning that the results in \eqref{Cor2_13} apply to symmetric RFC. As mentioned in \cite{BroDreKre_96} not all RFC methods are symmetric; however, for symmetric RFC the so-called Madelung rules apply, i.e., deletion pairs commute, meaning that it does not matter the order in which the sequences are deleted. However, if the primal concern is to apply this technique online, no deletion is actually possible since the estimation is done directly on measurements. 

\subsection{Crack Growth approaches}

Another alternative for fatigue estimation is the crack growth approach, which can be both addressed from a deterministic view-point using Paris' law (\cite{paris63}), or a stochastic perspective using for example jump processes, diffusion processes or stochastic differential equations (SDEs). However, in the crack growth approach a microscopic scale perspective is taken, thus making it difficult to transport to system level. We refer the interested readers to \cite{Sobczyk_FCG}, \cite{fatemi_98} and the references therein.

\section{Methods Comparison and Discussion} 

\begin{figure}[htb]
\centering
\includegraphics[scale=0.7]{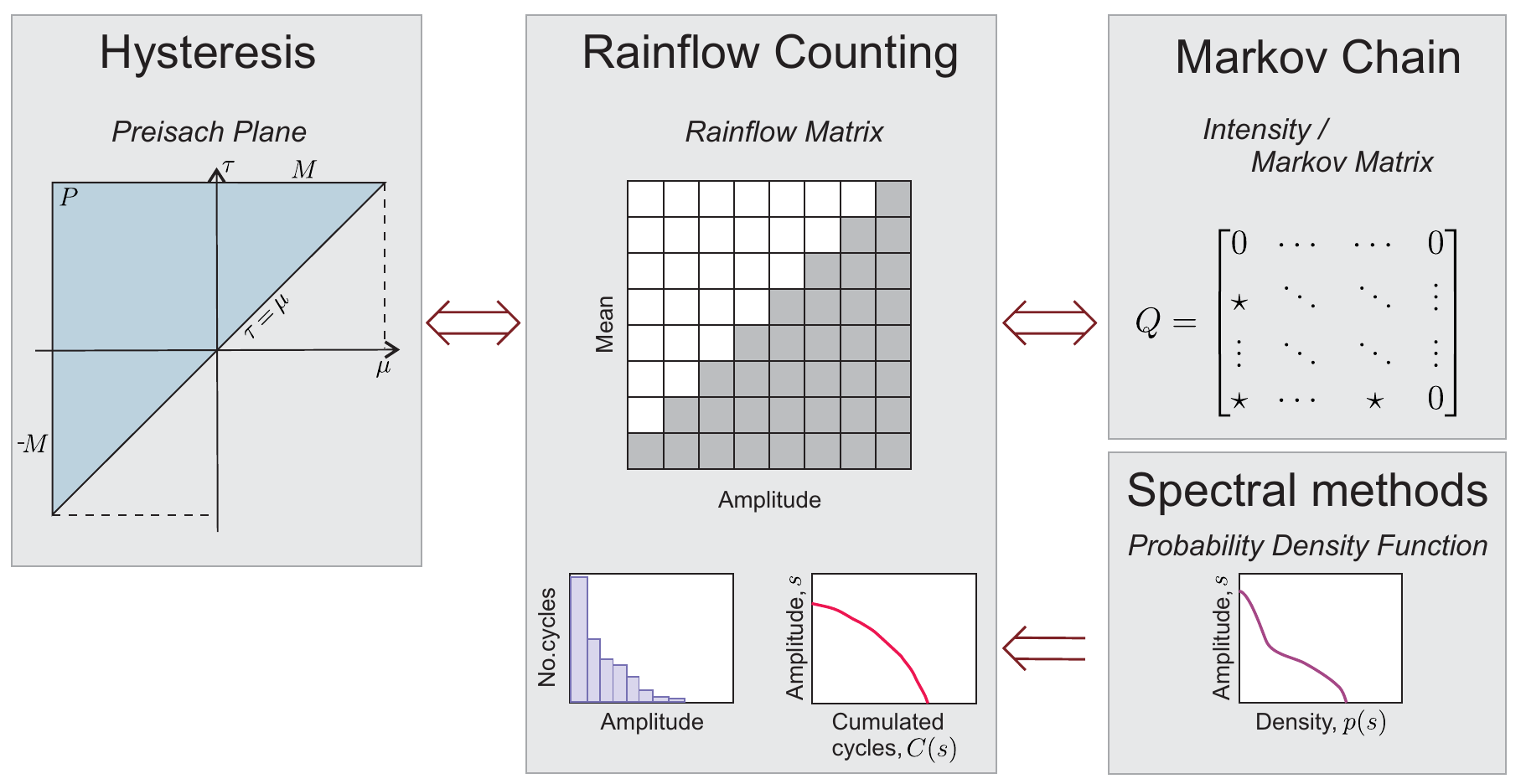}
\caption{Relationship between the compared methods.}
\label{fig:ConnectionDiagr}
\end{figure}

The aforementioned fatigue estimation methods share certain relations between each other. Firstly, there is an equivalence between the rainflow matrix and the Markov matrix or intensity of the Markov chain. Moreover, both have zeros below the diagonal, which is also the case for the Preisach plane $P$ in the Hysteresis method. The Spectral methods are related to RFC, since their intention is to approximate the rainflow density by spectral formulas, and they also relate to the stochastic methods in that their goal is to approximate certain density function. The hysteresis method is strongly related to the RFC, since the RFC actually identifies the closed hysteresis loops by counting cycles. A sketch of these relationships is depicted on Figure \ref{fig:ConnectionDiagr}.

Furthermore, a method comparison summary is shown on Table \ref{tab:AdvDisadv}, where advantages and disadvantages are presented for each method previously introduced.

\begin{table}[htb]
\begin{center} \begin{tabular}{l l l} 
\hspace{0pt}\textbf{Method} & \hspace{0pt}\textbf{Advantages}    &\hspace{0pt}\textbf{Disadvantages} \\ \hline
Rainflow          & Active Standard (ASTM E1049) & Post-processing \\      
\;Counting      	& Widely used                & Relies on linear accum. hypothesis \\
									&                            & Algorithmic, very non-linear \\ \hline
Spectral          & Can be used for control    & Black-box                    \\
									& Based on statistical measures & Narrow-band approximation \\ \hline
Stochastic        & Account for random loading & Parameters generally unknown \\ 
\;Methods 				& Could be used for prediction & May involve PDEs, SDEs       \\ 
									&                            & Very abstract formulation \\ \hline
Hysteresis        & Online estimation          & Typically hard control problem \\
									& Strong physical interpretation  & Density generally unknown \\
                  & Close mathematical form    & Approximation may be needed  \\ 
									&                            &                              \\ 
\end{tabular}
\caption{Methods advantages and disadvantages.}
\label{tab:AdvDisadv}
\end{center}\end{table}

\begin{figure}[htb]
\centering
\includegraphics[scale=0.8]{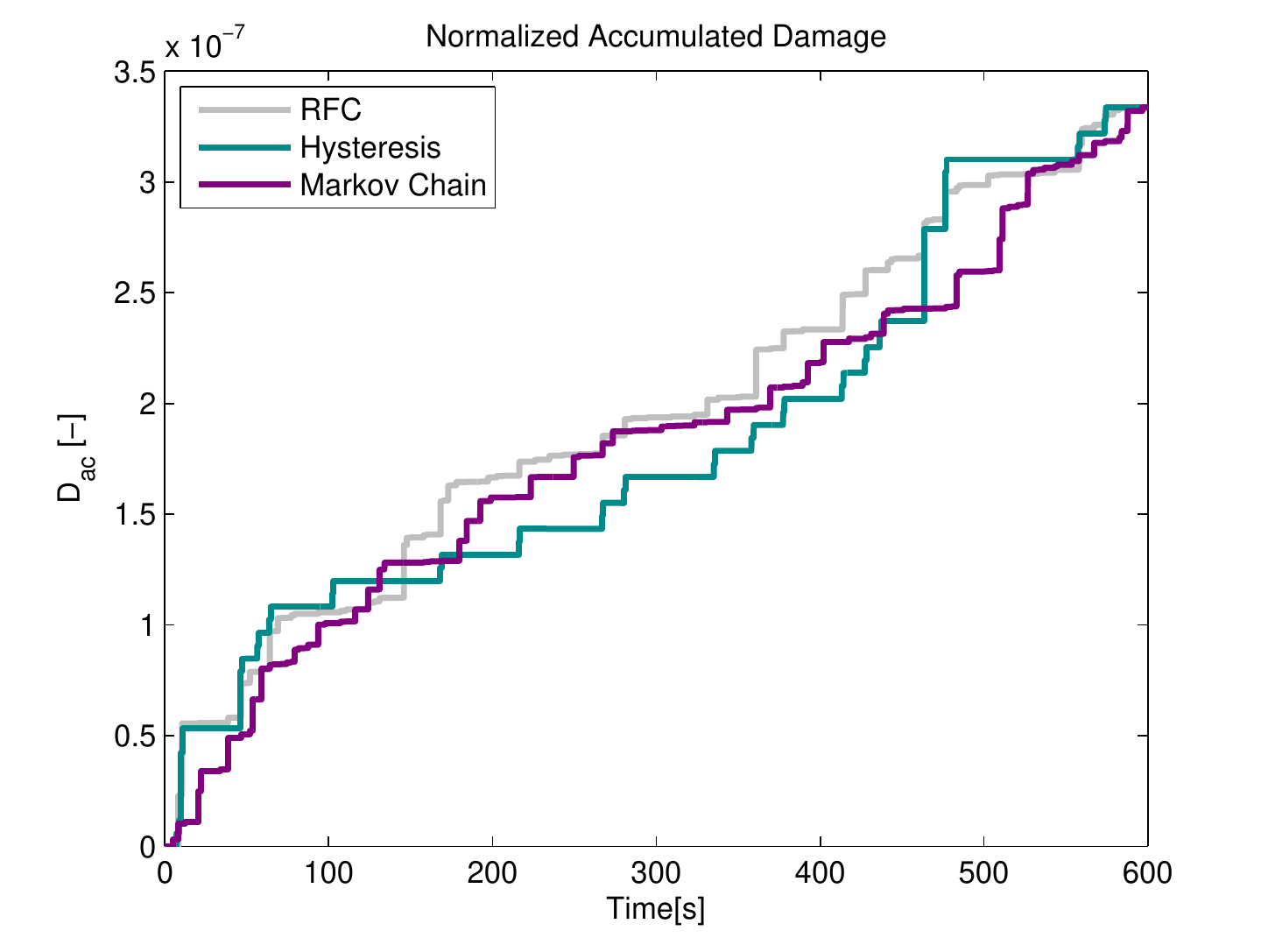}
\caption{Normalized accumulated damage for different estimation methods.}
\label{fig:RFCHystMC}
\end{figure}

For the next comparison part we will focus just on the MC instead of the whole stochastic methods class, which is quite broad. The accumulated damage provided by the RFC, MC and Hysteresis methods are compared in Figure \ref{fig:RFCHystMC}. The damage given by the hysteresis was normalized, such that it matches the accumulated damage of the RFC. The spectral method example could not be included, since the method delivers the damage rate itself and not instantaneous measurements. For the RFC and the MC method presented here, the instantaneous damage is given every time an extrema occurs and zero elsewhere, which is exactly what the hysteresis does, i.e., hold the value between certain thresholds. All these techniques can be used as post-processing tools, however not all of them can be used in the control loop. A brief summary is presented on Table \ref{tab:Control}, where it is reported if the methods can be implemented directly online or indirectly, i.e., not using measurements. The spectral methods are included indirectly, since they were included in the loop through transfer functions and not based on measurements. The Markov chain could be included online if the intensity matrix is parametrized with respect to the controls, which may not be realizable.

\begin{table}[htb]
\begin{center} \begin{tabular}{l c c l}
\textbf{Method} & \textbf{Online}  &\textbf{Indirect} &\textbf{Comments} \\ \hline
RFC               & - & -  & Only Post-processing   \\ \hline
Spectral          & - & X & Moments obtained by transfer function  \\ \hline
%Markov Chain      & X & X & Depends on intensity    \\ \hline
Hysteresis        & X & - & Approximation may be needed   \\ 
                  &   &   &                               \\ 
\end{tabular}
\caption{Methods applicability for control.}
\label{tab:Control}
\end{center}\end{table}

\section{Conclusions}

The literature regarding fatigue estimation methods is vast, since fatigue is an entire discipline by itself. The aim of the present paper is to provide a guide to the most recognized methods, which were assembled in four groups. These methods were presented and compared, from a control perspective in a Wind Turbine setting by estimating the damage from a tower bending moment time-series. A chart describing their advantages and disadvantages is presented on Table \ref{tab:AdvDisadv} and their applicability to control in Table \ref{tab:Control}. We also attempted to shed some light on the underlying relations between them.

Summarizing, the most widely used and standardized method is the RFC, but its algorithmic nature restricts its usage primarily as a post-processing tool. The spectral methods provide an alternative by trying to emulate the rainflow density function, they are based on statistical measures that are easier to calculate, but they are black-box and restricted (mainly) to narrow-band processes. The stochastic methods can accommodate the randomness of fatigue, but their construction is abstract and complicated, often involving stochastic or partial differential equations, and their parameters may need identification. The hysteresis method can be implemented online, acting on instantaneous measurements, but its complex and non-linear nature results in hard control problems. In general, one could say that the controller will influence the loading in the wind turbine components, and thus for implementing any of these techniques in the control loop, variable load should be considered by the estimation method in some sense.

\section{ACKNOWLEDGEMENT}
This work was partially supported by the Danish Council for Strategic Research (contract no. 11-116843) within the `Programme Sustainable Energy and Environment', under the ``EDGE'' (Efficient Distribution of Green Energy) research project.

%\ack This work was partially supported by the Danish Council for Strategic Research (contract no. 11-116843) within the `Programme Sustainable Energy and Environment', under the ``EDGE'' (Efficient Distribution of Green Energy) research project.

%\begin{thebibliography}{9}
%
%\bibitem{R1} Kopka~H, Daly~PW. 2003. \emph{A Guide to \LaTeX} (4th~edn).
%Addison-Wesley.
%
%\bibitem{R2} Lamport~L. 1994. \emph{\LaTeX: a Document Preparation System} (2nd~edn).
%Addison-Wesley.
%
%\bibitem{R3} Mittelbach~F, Goossens~M. 2004. \emph{The \LaTeX\ Companion}
%(2nd~edn). Addison-Wesley.
%\end{thebibliography}

\bibliographystyle{IEEEtran}
\bibliography{WindEnergyRefs}		

\end{document}